\input amssym.def
\input amssym.tex
\input graphicx.tex
\magnification=1100
\hsize=13.4cm \vsize=20.2cm
\parindent=0pt
\def\fips#1 {\midinsert\centerline{\includegraphics{fig#1.eps}}\medskip  
\centerline{\rmsmall Fig.\ #1}\endinsert}
\def\litem#1{\par \noindent\hangafter=0\hangindent=14pt
\llap{\hbox  to 14pt{#1\hfil}}\ignorespaces}
\font\sans=cmss10 at 7pt
\font\itsmall=cmmi8
\font\slsmall=cmsl8
\font\ttsmall=cmtt8
\font\rmsmall=cmr8
\font\title=cmr10 at 14pt
\def\tr{\raise4pt\hbox{\sans T}}
\def\bib #1#2#3#4{\litem{[#1]}{#2: {\slsmall #3} #4\par }}
\long\def\boxit#1{\par \setbox4=\hbox to12.8cm{\hsize=12.8cm\vbox{#1\par}}
\kern4pt\vbox{\hrule\hbox{\vrule\kern7pt
\vbox{\kern6pt\box4\kern6pt}
\kern6pt \vrule}\hrule}}
{\title A geometric approach to the diophantine Frobenius problem}
\vskip1cm Christian Blatter
 \vskip1cm\baselineskip=9.5pt
 {\rmsmall  ABSTRACT. It turns out that all instances of the diophantine Frobenius problem for three coprime {\itsmall a}$_i$ have a common geometric structure which is independent of arithmetic coincidences among the {\itsmall a}$_i$. By exploiting this structure we easily obtain   Johnson's formula for  the largest  non-representable {\itsmall z}, as well as  a formula for the number of such {\itsmall z}.  A procedure is described which computes these quantities in {\itsmall O}$($log(max\thinspace{\itsmall  a}$_i$)$)$ steps.}
 \vskip1cm
 \baselineskip=12pt
 1. INTRODUCTION \bigskip 
For an $n$-tuple ${\bf a}=(a_1, a_2, \ldots, a_n)$ of positive integers we denote by $T:=T({\bf a})$ the set of natural numbers $z$ that can be written in the form $$z=\sum_{k=1}^n x_k\,a_k\>,\qquad x_k\in\Bbb N:=\{0,1,2\ldots\}\>,$$
 and by $F:=F({\bf a})$  the set of natural numbers $z$ that cannot be so represented. If ${\rm gcd}(a_1,\ldots,a_n)=1$ then it is easily seen that all sufficiently large numbers $z$ are in $T$. It follows that in this case $F$ is a finite set, so there is a largest non-representable number $g({\bf a}):=\max F$. To compute this number and  maybe even the cardinality $N({\bf a})$ of  $F$  in terms of $a_1$, $\ldots$, $a_n$ constitutes the so-called {\sl diophantine Frobenius problem}.  We recommend [4], printed in 2005, as a comprehensive source of material about this problem; the  bibliography alone contains about  500 items.   \bigskip 
The case $n=2$ was first considered and solved by Sylvester [6], [7]. He proved: \medskip 
 {\bf Proposition 1.} {\sl If $a_1$, $a_2$ are $>1$ and coprime then}
 $$g(a_1,a_2)=a_1\,a_2 - a_1-a_2\>,\qquad N(a_1,a_2)={(a_1-1)(a_2-1)\over 2}\ .\eqno(1)$$
The proof of this result follows from inspection of Fig.~1 and is given at the beginning of the next section. 
 \bigskip 
 This paper deals with the case $n=3$. We shall give a natural geometric description of the set $F$ from which an explicit answer to the Frobenius-3-problem can immediately be read off. In order to formulate our result we introduce the quantities
 $$l_3:=\min\{l \in \Bbb N_{>0}\,|\, l\,a_3\in T(a_1,a_2)\}\qquad\circlearrowright\ .$$
 Here and in the sequel the sign $\circlearrowright$ indicates that there are three such formulae in all, whereby the other two are obtained by cyclic permutation $1\to2\to3\to1$ of the indices. About the $l_i$  the following can be said right away (cf.\ [2], Theorem 3 and eq.\ 26): \medskip 
 {\bf Proposition 2.}  {\sl Assume that the three numbers $a_1$, $a_2$, $a_3$ are pairwise prime and that $l_i\geq2$ $\>(1\leq i\leq3)$. Then the minimal representation 
  $$l_3\, a_3= x_{31}\, a_1 + x_{32}\,a_2\>, \quad x_{31}, x_{32}\in\Bbb N\qquad\circlearrowright \eqno(2)$$
   of  $\,l_i\, a_i$ $\,(1\leq i\leq 3)$ is uniquely determined, and one has
   $$x_{ij}\geq 1\qquad({\rm all}\ i\ne j)\ . \eqno(3)$$
Furthermore the $l_i$ are coupled to the $x_{ij}$ through}
  $$l_3=x_{13}+x_{23}\qquad\circlearrowright \ .\eqno(4)$$
  \medskip 
We now  state our main result; it will be proven in section 3:
  \medskip 
  {\bf Theorem 3.} {\sl Assume that the three numbers $a_1$, $a_2$, $a_3$ are pairwise prime and that $l_i\geq2$ $\>(1\leq i\leq3)$. Then}
  $$g({\bf a})=l_1\, l_2\, l_3 +\max \{x_{12} x_{23} x_{31},\> x_{21} x_{32} x_{13}\}-\sum\nolimits_i a_i\>;\eqno(5)$$
  $$N({\bf a})={1\over2}\biggl(\sum_i (l_i-1) a_i\> -\>l_1\,l_2\, l_3+1\biggr)\ .\eqno(6)$$ \medskip 
   The assumption $l_i\geq2$ means that none of the $a_i$ is ``superfluous''. If, e.g., $l_3=1$ then $F(a_1,a_2,a_3)= F(a_1,a_2)$; this case  is handled in Proposition 1. When the given $a_i$ are not pairwise prime then there is a way to get rid of  common factors, see [2], Theorem~2. Our formula (5), resp.\ its preliminary version (8), appears as Theorem 4 in [2] and on p.~35 of [4]. Note, however,  that  the proof given in [4] uses heavy algebraic machinery and is deferred to a later chapter.
 \medskip 
{\sl Example\/}\ 1. Let $a_1:=2n-1$, $a_2:=2n$, $a_3:=2n+1$ for an $n\geq2$. As $a_2\equiv 1$ and $a_3\equiv2$ (mod $a_1$) the smallest multiple of $a_1$ in $T(a_2,a_3)$ is $a_2+(n-1) a_3$;  similarly the smallest multiple of $a_3$ in $T(a_1,a_2)$ is $n a_1+a_2$, and obviously the minimal representation of $a_2$ is $2a_2=a_1+a_3$. Altogether we have
$$l_1=n+1, \> x_{12}=1, \> x_{13}=n-1; \> l_2=2, \>x_{21}=1, \>x_{23}=1; \> l_3=n, \> x_{31}=n, \> x_{32}=1\>;$$
so Theorem 3 gives
$$\displaylines{g({\bf a})=2n(n+1)+\max\{n,n-1\}- 6n=2n^2-3n\>,\cr
 N({\bf a})={1\over2}\bigl(n(2n-1)+2n+(n-1)(2n+1)-2n(n+1)+1\bigr)=n^2-n\ .\cr}$$
 \medskip 
For the computation of the $l_i$ and the $x_{ij}$ we propose the so called Lagrange algorithm -- a kind of   two-dimensional euclidean algorithm modeled after a Gram-Schmidt-process -- which takes $O\bigl(\log(\max a_i)\bigr)$ steps. The resulting procedure is developped in sections 4 and 5 of this paper. In [4] several other algorithms for $g({\bf a})$ are described, notably the algorithm of R\o dseth [5] which is an improved version of an earlier continued fraction algorithm by Selmer \& Beyer. \bigskip   \goodbreak
   2. PRELIMINARIES  \medskip   \fips 1  \medskip 
{\sl Proof of Proposition\/} 1.  We draw in the $(x_1,x_2)$-plane the directed graph $\Gamma$ with vertex set $\Bbb Z^2$ and edges of unit length connecting neighboring lattice points in the direction of increasing $x_1$ resp.\ $x_2$, see Fig.~1. For given $a_1$, $a_2\in\Bbb N_{>1}$ we define the {\sl height function} 
  $$f(x_1, x_2)\>:=\>a_1\, x_1+a_2\, x_2\ .$$
  Two  points in $\Bbb Z^2$ have the same height iff they differ by a vector ${\bf u}\in L:= \Bbb Z^2\cap f^{-1}(0)$. Since $a_1$ and $a_2$ are coprime  the set $L$ is the  one-dimensional lattice formed by the vectors ${\bf u}_k:=(k a_2, -ka_1)$, $\,k\in\Bbb Z$. \medskip 
  An integer $z\in\Bbb N$ can be represented in the form $z=x_1\, a_1 +x_2\,a_2$ with $x_1$, $x_2\in\Bbb N$ iff there is a directed edge path in $\Gamma$ connecting a point ${\bf u}_k\in L$   with a point ${\bf x}\in\Bbb Z^2$ of height $f({\bf x})=z$. Now the lattice points  that can be reached from a given ${\bf u}_k\in L$ lie  in the set $Q_k:={\bf u}_k+\Bbb R_{\geq0}^2\,$, a first quadrant with origin at ${\bf u}_k$, and  the set of all possible end points of such paths consists of the lattice points in  the union $\Omega:=\bigcup_k Q_k$ of these quadrants.  \medskip
    The lattice points of positive height that cannot be reached from one of the points ${\bf u}_k$ are the interior lattice points of the rectangular triangles $\Delta_k$ with vertices ${\bf u}_k$, ${\bf u}_{k+1}$, ${\bf u}_k+(a_2,0)$. The largest occurring height in such a $\Delta_k$  is given by the first formula (1), and using symmetry  one concludes that each $\Delta_k$ contains exactly $(a_1-1)(a_2-1)/2$ lattice points in its interior, which all have different heights. \hfill$\square$   \bigskip   
     For later purposes we note the following: Any lattice point in $\Delta_0$ can be connected by an admissible path to the point $(a_2,0)$ of height $a_1 a_2$. It follows that a number $z>0$ is in $F(a_1,a_2)$ iff there are integers $k_1$, $k_2\geq1$ such that $z=a_1a_2- k_1 a_1-k_2 a_2$.  \bigskip 
{\sl Proof of Proposition\/} 2.   First we show (3). Assume, e.g., that $x_{13}=0$. Then we have $l_1a_1=x_{12}a_2$, and as $a_1$, $a_2$ are coprime it follows that $l_1\geq a_2$. On the other hand, from $l_3\geq 2$, i.e., $a_3\in F(a_1,a_2)$ it follows that there are integers $k_1$, $k_2\geq1$ with $$a_3=a_1a_2-k_1 a_1-k_2 a_2 =(a_2-k_1)a_1-k_2 a_2\ .$$ Whence we would have $(a_2-k_1)a_1=k_2 a_2 +a_3$ which would imply $l_1<a_2$ -- a con\-tradiction. \medskip 
  We next show (4). Let the $x_{ij}$ be determined such that (2) holds. Then the quantities
  $$\mu_3:=l_3-x_{13}-x_{23}\qquad\circlearrowright$$
  satisfy $\mu_1a_1+\mu_2a_2+\mu_3a_3=0$. If the $\mu_i$ do not
  all vanish then up to a permutation of the  $a_i$ we must have one of the following:
  $${\rm (a)}\qquad \mu_1>0\>, \quad \mu_2<0\>,\quad \mu_3\leq0\>;\hskip1.5cm
  {\rm (b)}\qquad\mu_1>0\>, \quad \mu_2>0\>,\quad \mu_3<0\ .$$
  In case (a) it follows that $$(l_1-x_{21}-x_{31})a_1=\mu_1a_1=(-\mu_2)a_2+(-\mu_3)a_3\>,$$ contradicting the definition of $l_1$. In case (b), from
  $$(x_{13}+x_{23}-l_3)a_3=-\mu_3 a_3=\mu_1a_1+\mu_2 a_2$$
  it follows by definition of $l_3$ that $x_{13}+x_{23}\geq 2 l_3$, whence, e.g., $x_{13}\geq l_3$. By definition of the $x_{ij}$ we now  have the representation
  $$(l_1-x_{31})a_1=(x_{12}+x_{32})a_2+(x_{13}-l_3)a_3$$
  which again contradicts the definition of $l_1$.  \medskip 
 As (4) is true for all possible choices of $x_{21}$ and $x_{31}$ consistent with the definition of $l_2$ and $l_3$, and as these choices can be made independently, it follows that there is in fact no choice at all, which means that the $x_{ij}$ are indeed uniquely determined. \hfill$\square$ \bigskip  \medskip 
 3. PROOF OF THE MAIN RESULT \bigskip
 We now come to the proof of Theorem 3. Inspired by the proof of Sylvester's result for $n=2$ we  embed the problem into the following geometric setup: Consider the integer lattice $\Bbb Z^3$ in  euclidean $(x_1,x_2,x_3)$-space $\Bbb R^3$. We use $\Bbb Z^3$ as vertex set of a directed graph $\Gamma$ whose edges are the segments of unit length connecting neighboring lattice points in the direction of increasing $x_1$, resp.\ $x_2$, $x_3$. Given $a_1$, $a_2$, $a_3$, we again define  the {\sl height function\/}
  $$f(x_1,x_2,x_3):=a_1\, x_1+a_2\, x_2 + a_3\, x_3$$
  which on the one hand is just a linear functional on $\Bbb R^3$ and on the other hand assigns a height $f({\bf x})$ to each lattice point ${\bf x}\in\Bbb Z^3$. The kernel $H:=f^{-1}(0)$ of $f$ is a plane through the origin of $\Bbb R^3$ and contains the {\sl Frobenius lattice\/} $L:=H\cap\Bbb Z^3$ of integer solutions to the equation $f({\bf x})=0$.  \medskip 
{\bf Lemma 4.} (a) {\sl Let $m_1a_1+m_2 a_2=1$ with $m_i\in\Bbb Z$. Then the vectors ${\bf e}_1:=(a_2, - a_1, 0)$, ${\bf e}_2:=(a_3 m_1,a_3m_2, -1)$ form a basis of $L$.}  \smallskip 
(b) {\sl  The fundamental domain of the  lattice $L$, when projected to  the plane $x_i=0$, has area $a_i$ $\>(1\leq i\leq 3)$.} \medskip 
{\sl Proof.} (a) One easily checks that $f({\bf e}_1)=f({\bf e}_2)=0$, which means that ${\bf e}_1$, ${\bf e}_2\in L$. On the other hand, let ${\bf u}=(u_1,u_2,u_3)$ be an arbitrary point of $L$. Then ${\bf u}+ u_3{\bf e_2}=(u_1',u_2',0)\in L$ which implies $a_1 u_1'+a_2 u_2'=0$. Since $a_1$, $a_2$ are coprime it follows that $(u_1',u_2',0)=k{\bf e}_1$ for a $k\in\Bbb Z$, whence we have ${\bf u}=k{\bf e}_1 -u_3{\bf e}_2$.  \smallskip 
(b) It suffices to compute the vector product
$${\bf e}_1\times {\bf e_2}=\bigl(a_1,a_2,(m_1 a_1+m_2 a_2) a_3\bigr)=(a_1,a_2,a_3)\ .\eqno\square$$
 \smallskip  
The three lattice  vectors
$${\bf f}_1:=(l_1,-x_{12}, -x_{13}), \ {\bf f}_2:=(-x_{21}, l_2, -x_{23}), \ {\bf f}_3:=(-x_{31}, -x_{32}, l_3)\>\in L$$
encoding the data $l_i$, $x_{ij}$ will play a special r\^ole. We shall call any  vector of the form  ${\bf f}_i$ or $-{\bf f}_i$ a {\sl basic vector\/} and any set of three essentially different vectors among  the  $\pm\,{\bf f}_i$ a {\sl solution basis\/} for the Frobenius problem at hand.
 
 \fips 2  \medskip 

  The connection of the three-dimensional structure described so far with the diophantine Frobenius problem is the following: A natural number $z$ is in $T({\bf a})$ iff there is a lattice point ${\bf u}\in L$ and a directed edge path in $\Gamma$ beginning at ${\bf u}$ and ending in a point ${\bf x}\in \Bbb Z^3$ of height $f({\bf x})=z$. Now the lattice points  that can be reached from a given ${\bf u}\in L$ lie  in the set $O_{\bf u}:={\bf u}+\Bbb R_{\geq0}^3\,$, an octant with origin at ${\bf u}$, and the set of  all possible end points of such paths consists of the lattice points in  the union $\Omega:=\bigcup_{\bf u\in L} O_{\bf u}$ of these octants. It follows that $T({\bf a})=f(\Bbb Z^3\cap \Omega)$. The   boundary    $\partial\Omega$ is  $L$-periodic; it consists of three L-shapes per octant and touches  the plane  $H$ in the points of $L$. It looks like a washboard and  is depicted in Fig.~2, {\sl as seen from below.} The proof of Theorem 3 essentially consists in understanding this figure. \medskip 
  Consider, e.g., the positive $x_3$-axis. It is an edge of the octant  $O_{\bf 0}$ and belongs to the boundary of $\Omega$ until it is intercepted at ${\bf p}=(0,0, u_3)$ by a wall $x_3=\>$const.\ belonging to another octant $O_{\bf u}$, ${\bf u}=(u_1,u_2,u_3)\in L$. For this to happen it is necessary that $u_1\leq 0$, $u_2\leq0$, $u_3>0$, and that there is no point ${\bf u}'=(u_1', u_2', u_3')\in L$ with $u_1'\leq0$, $u_2'\leq0$ and $u_3'<u_3$. This means
  $$\eqalign{u_3&=\min\bigl\{u\in \Bbb N_{>0}\,\bigm|\,\exists\, u_1, u_2\in \Bbb Z_{\leq0}\!: \>  a_1 u_1+a_2 u_2 + a_3 u=0\bigr\}\cr &=\min\bigl\{u\in \Bbb N_{>0}\,\bigm|\,u\, a_3\in T(a_1, a_2)\bigr\}\>,\cr}\eqno(7)$$
  from which we deduce
  $$u_3=l_3\>, \quad -u_1=x_{31}\>, \quad -u_2=x_{32}\>,$$
  i.e., $(u_1, u_2, u_3)={\bf f}_3\,$; and similarly for the other $l_i$, $x_{ij}$. \medskip 
  The lattice points ${\bf x}$ of positive height that cannot be reached from a point ${\bf u}\in L$ are the interior lattice points contained in the region $W$ enclosed between $H$ and $\partial\Omega$. As seen in the figure, the restriction $f\!\restriction\!\partial\Omega$ takes local maxima at the points  ${\bf q}_1$, ${\bf q}_2$, and the interior lattice points of maximal height are  ${\bf q}_1-(1,1,1)$ or ${\bf q}_2-(1,1,1)$ and their equivalents mod$\,L$. It follows that the maximal non-realizable height $g({\bf a})$ is given by
  $$g({\bf a})= \max\{f({\bf q}_1), f({\bf q}_2)\}-\sum\nolimits_i a_i 
  =l_3\,a_3+\max\{x_{21}a_1, x_{12}a_2\} -\sum\nolimits_i a_i\ .\eqno(8)$$  \medskip 
  Since the three L-shapes have areas $a_i$ by Lemma 4(b), we deduce from Fig.~2 that the $a_i$ satisfy
  $$a_1=x_{12}l_3 + x_{13}x_{32}\>, \quad a_2=x_{21}l_3+x_{23}x_{31}\>, \quad a_3=l_1 l_2-x_{12}x_{21}\qquad\circlearrowright\ .\eqno(9)$$
Substituting these expressions into (8) one arrives at the symmetric formula (5):
  $$\eqalign{g({\bf a}) &=
  l_1 l_2 l_3- l_3 x_{12}x_{21}+\max\{x_{21}x_{12} l_3+x_{21}x_{13}x_{32},
  x_{12}x_{21}l_3+x_{12}x_{23}x_{31}\}-\sum\nolimits_i a_i \cr
  &=l_1 l_2 l_3 +\max\{x_{21}x_{13}x_{32},
  x_{12}x_{23}x_{31}\}-\sum\nolimits_i a_i\ . \cr}$$
 \smallskip   We now come to the proof of formula (6). We have to count the number $z_0$ of interior lattice points in the quotient $\hat W:=W/L$. As $\hat W$ does not have a simple description in terms of inequalities we are going to determine $z_0$ ``from the outside'' by means of a three-dimensional analog of Pick's area formula. Let $z_1$ denote the total number of relative interior lattice points in the three L-shapes; similarly, let $z_2$ be the total number of relative interior lattice points on the reentrant edges of $\hat W$ and $z_3$ be the number of such points on the protruding edges of $\hat W$. Then we have the following formula:  
  $${\rm vol}(\hat W)=z_0+{1\over2} z_1 +{3\over4} z_2 +{1\over 4} z_3+{7\over4}\>,
  \leqno{\bf Lemma\ 5.}$$   
  whereby the last term incorporates the contribution of the corners of $\hat W$. \medskip 
 {\sl Proof.} We perform a ``Gedankenexperiment'' used already in [1] for a proof of Pick's area formula. Imagine that at time $0$ a unit of heat is concentrated at each point of $\Bbb Z^3$.   This heat will be distributed all over space by heat conduction, and at time $\infty$ it will be equally distributed in space with density $1$. In particular, the amount of heat contained in $\hat W$ will be ${\rm vol}(\hat W)$. Where does this amount of heat come from? For symmetry reasons there is absolutely no flux across the unit squares of $\partial\Omega$, and, again by symmetry, the net flux across $H/L$ is $0$ as well. As a consequence, the final amount of heat within $\hat W$ comes from the interior lattice points, counted by $z_0$, and from the lattice points on the boundary of $\hat W$. The lattice points counted by $z_1$ send half their heat into $\hat W$, whereas the corresponding factor is ${3\over4}$ for the points counted by $z_2$ and ${1\over4}$ for the points counted by $z_3$. Furthermore $\hat W$ possesses two protruding corners ${\bf q}_1$ and ${\bf q}_2$ which contribute ${1\over8}$ each, 
 three ``L-corners'' contributing ${3\over8}$, and finally the reentrant corner on $H$ which contributes ${3\over8}$ as well.\hfill$\square$ \medskip 
  \fips 3  \medskip 
  In order to compute ${\rm vol}(\hat W)$ directly we use the L-shape $A$ in the plane $x_3=0$ as  fundamental domain. Fig.~3 shows  $A$, as seen from the positive $x_3$-direction. The vertical prism $K$ determined by $A$ and the plane
  $$H:\qquad x_3=-{1\over a_3}(a_1 x_1 +a_2 x_2)$$
is a representative for $\hat W$. To compute the volume of $K$ we split $A$ into two rectangles  and work with the heights of $K$ in their midpoints.  We obtain
  $${\rm vol}(\hat W)=
  {x_{31}x_{12}\over a_3}\Bigl(a_1{x_{31}\over2}+a_2{2x_{32}+x_{12} \over2}\Bigr) +
  {l_1 x_{32} \over a_3}\Bigl(a_1{l_1 \over 2} + a_2{x_{32}\over2}\Bigr)\>,$$
  which using (2) and (9) can be brought into the following symmetric form:
  $${\rm vol}(\hat W)={1\over2}\left(\sum\nolimits_i l_i\,a_i\>- l_1l_2l_3\right)\ .\eqno(10)$$
  We now compute the quantities $z_1$, $z_2$ and $z_3$. -- The horizontal L-shape $A$ in Fig.~3 has area $a_3$ by Lemma 4(b)  and $2(l_1+l_2)$ lattice points on its boundary. Therefore by Pick's area formula for the plane the number of interior lattice points on $A$ is given by $a_3-(l_1+l_2)+1$, and we obtain
  $$z_1=\sum\nolimits_i a_i - 2\sum\nolimits_i l_i+3\ .\eqno(11)$$
  Inspection of Fig.~2 shows that
  $$z_2=\sum\nolimits_i(l_i-1)=\sum\nolimits_i l_i -3\>, \qquad z_3=\sum\nolimits_{i\ne j} (x_{ij}-1)=\sum\nolimits_i l_i -6\ .\eqno(12)$$
  Introducing (10), (11) and (12) into Lemma 5 we get
  $$\eqalign{z_0&={\rm vol}(\hat W)-{1\over2} z_1 -{3\over4} z_2 -{1\over 4} z_3-{7\over4} \cr
  &={1\over2}\left(\sum\nolimits_i l_i\,a_i\>- l_1l_2l_3\right)-{1\over2}\left(\sum\nolimits_i a_i - 2\sum\nolimits_i l_i+3\right)-{3\over4}\left(\sum\nolimits_i l_i -3\right)\cr &\quad -{1\over 4} \left(\sum\nolimits_i l_i -6\right)-{7\over4} \cr
  &= {1\over2}\left(\sum\nolimits_i (l_i-1)a_i-l_1 l_2 l_3+1\right)\ .\cr}$$ 
  This concludes the proof of Theorem 3. \hfill$\square$\bigskip  \medskip 
  4. FINDING A SOLUTION BASIS  \bigskip 
In order to make Theorem 3 useful we have to establish a procedure to  compute the quantities $l_i$, $x_{ij}$. Our argument takes place in the plane $H$. Fig.~4 shows $H$ as seen from the tip of the vector ${\bf a}$, the points of $L$ are again marked by bullets. The three planes $x_i=0$ intersect $H$ in three lines $g_i$ through the origin which altogether divide $H$ into six sectors of various widths. The line $g_3$, spanned by the vector ${\bf v}_3:=(-a_2, a_1,0)$, is at the same time a level line of the linear function $x_3$ restricted to $H$. Equation (7) can now be interpreted as follows: In order to find $l_3$ we have to translate the line $g_3\!: \,x_3=0$ in the direction of increasing $x_3$ (marked by an arrow in Fig.~4)    until it hits for the first time a  lattice point in the sector $x_1\!\!\leq\!0\>\,\wedge\> x_2\!\leq\!0$. The lattice point obtained in this way is the point ${\bf f}_3$. Translating similarly the lines $g_1$  and $g_2$ one obtains the lattice points ${\bf f}_1$, ${\bf f}_2$ in the appropriate sectors.  \medskip 
  \fips 4  \medskip 
Being minimizers of some sort the ${\bf f}_i\in L$ tend to be short. Now there exists a revered algorithm (attributed to Lagrange, Gauss and others, see [3]) which finds the shortest vector ${\bf u}$ of the lattice $L$, and we plan to make use of this algorithm.  But if ${\bf u}$ happens to lie in a sector of width $>60^\circ$, as in Fig.~4, there is no guarantee that ${\bf u}$ coincides with the basic vector ${\bf f}_i$ (or $-{\bf f}_i$)  in that sector. For this reason we change the metric in such a way that the three lines $x_i=0$ intersect at angles of $60^\circ$.    \bigskip 
{\bf Lemma 6.} {\sl For a suitable  scalar product $\langle{\bf x},{\bf y}\rangle:={\bf x}\tr\, Q\,{\bf y}$
the three lines $x_i=0$ in $H$ intersect at angles of $60^\circ$.} \bigskip 
{\sl Proof.} The three directions in question are
$${\bf v}_1:=(0, -a_3, a_2)\>,\quad {\bf v}_2:=(a_3, 0, -a_1)\>,\quad {\bf v}_3:=(-a_2, a_1, 0)\ .$$
Consider now the linear map $P:\Bbb R^3\to \Bbb R^2$ given by the matrix
$$P:=\left[\matrix{-a_1 & a_2 &0 \cr 0 &0 & \sqrt{3} a_3 \cr}\right]\ .$$
The kernel of $P$ is spanned by the vector $(a_2,a_1,0)\notin H$, therefore the restriction $P\restriction H$ maps $H$ bijectively onto the euclidean plane $E:=\Bbb R^2$. One easily computes
$$P{\bf v}_1=2a_2a_3\Bigl(-{1\over2},{\sqrt{3}\over2}\Bigr)\>,\quad 
P{\bf v}_2=2a_1a_3\Bigl(-{1\over2},-{\sqrt{3}\over2}\Bigr)\>,\quad 
P{\bf v}_3=2a_1a_2\>(1,0)\>,$$ 
which shows that in the image plane the  lines $g_i$ intersect at  angles of $60^\circ$.  \smallskip  
Pulling back the euclidean scalar product in $E$ to $H$ one obtains there the new scalar product
$$\langle{\bf x},{\bf y}\rangle:=(P{\bf x})\tr\>P{\bf y}={\bf x}\tr P\tr\,P \>{\bf y}={\bf x}\tr\, Q \>{\bf y}\>,$$
where the integer matrix $Q:=P\tr\, P$ is given by
$$Q=\left[\matrix{a_1^2 &-a_1a_2 & 0 \cr - a_1 a_2 &  a_2^2 &0 \cr 0 & 0 &3a_3^2 \cr}\right]\ .$$
\vskip-12pt \rightline{$\square$}  \bigskip  
In what follows, $|\cdot|$ denotes the norm corresponding to the scalar product $\langle\cdot,\cdot\rangle$. \medskip 
Lagrange's algorithm (to be described in the next section) produces a {\sl reduced basis\/} $({\bf u}, {\bf v})$ of the Frobenius lattice $L$. This means that ${\bf u}$ is a shortest nonzero vector in $L$ and that ${\bf v}$ is a shortest vector in $L\setminus \Bbb Z\,{\bf u}$; in particular, $|{\bf u}|\leq |{\bf v}|$. \bigskip 
{\bf Theorem 7.} {\sl Under the hypotheses of Theorem\/} 3, {\sl let $({\bf u},{\bf v})$ be a reduced basis of the Frobenius lattice $L$. Assume that $u_1\leq0$, $u_2\leq 0$, $u_3>0$ and put $\lambda:=v_3/u_3$. Then the three vectors $${\bf u}\>, \quad {\bf v}_-:=
{\bf v}-\lceil\lambda\rceil{\bf u}\>, \quad {\bf v}_+:={\bf v}-\lfloor\lambda\rfloor{\bf u}$$   form a solution basis for the Frobenius problem determined by the $a_i$.} \medskip 
\fips 5  \medskip 
{\sl Proof.} We argue in the $(x,y)$-plane $E$, but omit the $P$ in our notation. The lines $g_i$ enclose angles of $60^\circ$, creating (closed) sectors $S_k$ $(1\leq k\leq 6)$, see Fig.~5.  After scaling we have ${\bf u}=(-\sin\theta,\cos\theta)$ with $|\theta|\leq{\pi\over6}$, whence ${\bf u}\in S_2$, and we may assume that the lattice line $\ell\parallel{\bf u}$  containing ${\bf v}$ is to the right of  ${\bf u}$. We shall show that the   points ${\bf u}$, ${\bf v}_-$ and ${\bf v}_+$  are basic points in the sectors $S_6$, $S_1$ and $S_2$ respectively. 
 \medskip 
  We begin with the remark that in fact $|\theta|<\pi/ 6$. Assume to the contrary that, e.g., $\theta=\pi/6$. Then ${\bf u}=(-1/2, \sqrt{3}/2)\in g_1$, meaning $u_1=0$. In this case ${\bf u}$ could not be basic by (3). It follows that there would have to be a lattice point in $S_2$ with  $y$-coordinate $<\sqrt{3}/2$. But there is no room for such a point since the interior of the unit circle is forbidden.
   \medskip 
 Let ${\bf b}:={\bf v}-\lambda{\bf u}=(b,0)$ be the point where $\ell$ intersects $g_3$. Then 
 $b=h/\cos\theta$ where $h$ denotes the distance from ${\bf 0}$ to $\ell$, whence $h\geq\sqrt{3}/2$.  We write $(\cos\phi,\sin\phi)=:{\bf e}_\phi$.   \bigskip 
{\bf Lemma\ 8.} \ (a)\quad 
$ \langle{\bf e}_{\pi/6},{\bf u}\rangle< \langle{\bf e}_{\pi/6}, {\bf b}\rangle\>,$\quad (b)\quad 
$\langle{\bf e}_{-\pi/6},-{\bf u}\rangle<
 \langle{\bf e}_{-\pi/6}, {\bf b}\rangle\ .$  \bigskip 
 {\sl Proof.} The left sides of (a) and (b) are $\sin\bigl({\pi\over 6}-\theta\bigr)$ and 
 $\sin\bigl({\pi\over 6}+\theta\bigr)$ respectively, so they both are $\leq \sin\bigl({\pi\over 6}+|\theta|\bigr)$. On the other hand the right sides of (a) and (b) both have the same value
 $$\cos{\pi\over 6}\> b=\cos{\pi\over 6}\>{h\over\cos\theta}\geq {3\over4\cos\theta}\ .$$
 It remains to prove that for $0\leq\theta<{\pi\over6}$ one has
 $$2 \sin\Bigl({\pi\over 6}+\theta\Bigr)\cos\theta< {3\over2}\ .$$
 But here the left side can be written as $\sin\bigl({\pi\over 6}+2\theta\bigr) +
 \sin\bigl({\pi\over 6}\bigr)$ which is $< 1 +{1\over2}$.\hfill$\square$
  \bigskip 
  Put ${\bf c}_-:={\bf b}-{\bf u}$, ${\bf c}_+:={\bf b}+{\bf u}$. Then from Lemma 8(a) it follows that $\langle{\bf e}_{\pi/6},{\bf c}_-\rangle>0$, whence ${\bf c}_-\in {\rm int}(S_6)$, and analogously from Lemma 8(b) it follows that $\langle{\bf e}_{-\pi/6},{\bf c}_+\rangle>0$, whence ${\bf c}_+\in {\rm int}(S_1)$.  From this we conclude that the lattice points ${\bf v}_-\in[{\bf b}, {\bf c}_-]$ and ${\bf v}_+ \in [{\bf b}, {\bf c}_+]$ lie in  $S_6$ and $S_1$ respectively.
   \medskip 
 As ${\bf c}_+\in  {\rm int}(S_1)$ the line $\ell$ intersects $g_2$ at a $y$-level $>\cos\theta$, whence all lattice points on $\ell\cap S_2$ have a larger $y$-level than ${\bf u}$, and similarly, as $-{\bf c}_-\in  {\rm int}(S_3)$, the line $\ell'$ intersects $g_1$ at a $y$-level $>\cos\theta$, whence all lattice points on $\ell'\cap S_2$ have a larger $y$-level than ${\bf u}$. This implies that the vector ${\bf u}$ is basic in  its sector $S_2$.    \medskip 
 Note that $\langle{\bf e}_{\pi/6},{\bf u}\rangle>0$, whence going upwards along $\ell\cap S_1$  the distance to $g_1$ increases. Since ${\bf v}_+$ is the first lattice point met along this path, ${\bf v}_+$ is basic for the sector $S_1$, unless there were an even better lattice point on the  parallel to $\ell$ through the point $2{\bf b}$. But the latter is prohibited by the inequality $\langle{\bf e}_{\pi/6},{\bf c}_+\rangle<\langle{\bf e}_{\pi/6}, 2{\bf b}\rangle$ which follows easily from Lemma 8(a). \medskip 
 Similarly one has $\langle{\bf e}_{-\pi/6},{\bf u}\rangle<0$, and this implies that going downwards along $\ell\cap S_6$ the distance to $g_2$ increases. Since ${\bf v}_-$ is the first lattice point met along this path, ${\bf v}_-$ is basic for the sector $S_6$, unless there were an even better lattice point on the  parallel to $\ell$ through the point $2{\bf b}$. But the latter is prohibited by the inequality $\langle{\bf e}_{-\pi/6},{\bf c}_-\rangle<\langle{\bf e}_{-\pi/6}, 2{\bf b}\rangle$ which follows easily from Lemma 8(b). \hfill$\square$
  \bigskip  \medskip 
  5. LAGRANGE'S ALGORITHM  \bigskip 
Lagrange's algorithm, as it is called in [3],   takes an arbitrary basis $({\bf u},{\bf v})$ of the Frobenius lattice $L$ as input and in a certain number of steps arrives at a reduced basis of $L$. An essential accessory to the calculations is the {\sl Gram matrix}
$$G:=G({\bf u},{\bf v}):=\left[\matrix{\langle{\bf u},{\bf u}\rangle & \langle{\bf u},{\bf v}\rangle \cr \langle{\bf v},{\bf u}\rangle &\langle{\bf v},{\bf v}\rangle \cr}\right]$$
of the current basis $({\bf u},{\bf v})$; it is updated along with the basis vectors. \medskip 
The following box is taken from [3]. The subscript $\leq$ to a basis indicates that  one assumes $|{\bf u}|\leq|{\bf v}|$, and $\lfloor\>\cdot\>\rceil$ denotes the nearest integer function. \medskip 
\boxit
{{\bf Input:} A basis $({\bf u}, {\bf v})_\leq$ with its Gram matrix $G=(g_{i,j})_{1\leq i,j\leq2}$. \par  
{\bf Output:} A reduced basis of $L$ with its Gram matrix. \smallskip 
\litem{1.} Repeat  \smallskip  
\litem{2.} \quad ${\bf r}:={\bf v}-x\,{\bf u}$, where $x:=\bigl\lfloor\, g_{1,2} / g_{1,1}\bigr\rceil\,$. 
\litem{} \quad
When computing $x$, also compute the remainder $y:= g_{1,2}- xg_{1,1}$.   \par  
\litem{3.} \quad ${\bf v}:={\bf u}\,$.  \smallskip 
\litem{4.} \quad ${\bf u}:={\bf r}\,.$  \smallskip 
\litem{5.} \quad Update the Gram matrix as follows: swap $g_{2,2}$ and $g_{1,1}$; then let $g_{1,2}:=y$ and \phantom{\quad}$g_{1,1}:=g_{1,1}-x(y+g_{1,2})$. \smallskip 
\litem{6.} Until $|{\bf u}|\geq|{\bf v}|\,.$  \smallskip 
\litem{7.} Return $({\bf v}, {\bf u})_\leq$ and its Gram matrix (setting $g_{2,1}:=g_{1,2}$).\par  } 
\centerline{\rmsmall Lagrange's algorithm}\bigskip 
We now combine this with the results of the foregoing section in order to obtain a coherent description of the computational procedure to determine the $l_i$, $x_{ij}$. \medskip 
 When $a_1$, $a_2$, $a_3$ are given, one first has to set up the basis $({\bf e}_1, {\bf e}_2)$ of $L$  given in Lemma 4(a). This requires $O\bigl(\log(\max_i a_i)\bigr)$ steps for the euclidean algorithm to find $m_1$, $m_2$. Using this basis as input one then runs Lagrange's algorithm and obtains a reduced basis $({\bf u}, {\bf v})_\leq$ of $L$. As shown in [3], Theorem 3.0.3, this is accomplished in at most $O\bigl(\log(\max_i a_i)\bigr)$ loops of the algorithm. Replacing ${\bf u}$ by $-{\bf u}$, if necessary, makes $u_i>0$, $u_j<0$ $(j\ne i)$ for some $i$. Now put $\lambda:=v_i/u_i$ and define ${\bf v}_-$, ${\bf v}_+$ as given in Theorem 7. The quantities $l_i$, $x_{ij}$ can then  be read off from the coordinates of the three vectors ${\bf u}$, ${\bf v}_-$, ${\bf v}_+$. Note however that the bit complexity of the whole computation is quadratic insofar as the bit-length of the input data $a_i$ not only affects the number of required steps/loops  but also the cost of each step. \bigskip 
{\sl Example\/}\ 2. Consider the random numerical example $a_1:=4327$, $a_2:=6716$, $a_3:=9237$.  In 8 steps the euclidean algorithm finds $2055 a_1- 1324 a_2=1$, and after 6 loops of Lagrange's algorithm we arrive at a reduced basis of $L$ given by ${\bf u}=(-53,-47,59)$, ${\bf v}=(-130,59,18)$. Theorem 7 then tells us that
$${\bf f}_1=(130, -59,-18)\>,\quad {\bf f}_2=(-77,106,-41)\>,\quad {\bf f}_3=(-53,-47,59)$$
is a solution basis for the given data, and by Theorem 3 we have
$$g({\bf a})=920\,947\>, \qquad N({\bf a})=493\,045\ .$$
 \bigskip {\bf References}\bigskip 
{\rmsmall \baselineskip=9.5pt\parskip=1pt
\bib{1}{Chr.\ Blatter} {Another proof of Pick's area theorem.}{Math.\ Mag.\ 70 (1997), 200.}
\bib{2}{S.M.\ Johnson}{A linear diophantine problem.}{Can.\ J.\ Math.\ 12 (1960), 390--398.}
\bib{3}{P.Q.\ Nguyen \& D.\ Stehl\'e}{Low-dimensional lattice basis reduction revisited.} {ACM Transactions on Algorithms 5 (2009), 4 (Oct.), 1--48.}
\bib{4}{J.L.\ Ramirez Alfonsin}{The diophantine Frobenius problem.} {Oxford University Press 2005}
\bib{5}{\O.J. R\o dseth}{On a linear diophantine problem of Frobenius.}{J.\ Reine und Angewandte Mathematik 301 (1978), 171--178.}
\bib{6}{J.J.\ Sylvester}{On subinvariants, i.e.\ semi-invariants to binary quantities of an unlimited order.}{Am.\ J.\ Math.\ 5 (1882), 119--136.}
\bib{7}{J.J.\ Sylvester}{Problem 7382.}{Educational Times 37 (1884), 26.} \bigskip 
\parindent=0pt \par 
{\slsmall Address of the author:} \smallskip 
Christian Blatter \par 
Department of Mathematics \par 
Swiss Fed.\ Inst.\ of Technology  \par 
CH-8092 Zurich (Switzerland)  \par  
e-mail: {\ttsmall christian.blatter@math.ethz.ch}}

   \end